\documentclass[reqno,12pt]{amsart}
\usepackage{hyperref}
\usepackage{amsmath}
\usepackage{amssymb} %
\usepackage{color}
\usepackage{xcolor}
\usepackage{ifpdf}
\usepackage{txfonts}
\usepackage{animate}
\usepackage{hyperref,url}
\usepackage{fleqn}
\usepackage[T1]{fontenc}
\usepackage[all]{xy}

\begin{document}

\title[\textbf{Non-isomorphic Pure Galois-Eisenstein Rings}]
{\textbf{Non-isomorphic Pure Galois-Eisenstein Rings}}

\author[Alexandre Fotue Tabue and Christophe Mouaha]
{Alexandre Fotue Tabue and Christophe Mouaha}

\address{Alexandre Fotue Tabue\newline
Department of mathematics, Faculty of Sciences,  University of
Yaounde 1, Cameroon} \email{alexfotue@gmail.com}

\address{Christophe Mouaha   \newline
Department of mathematics,  Higher Teachers Training College of
Yaounde, University of Yaounde 1, Cameroon}
\email{cmouaha@yahoo.fr}

\maketitle \numberwithin{equation}{section}
\newtheorem{Theorem}{Theorem}[section]
\newtheorem{Example}{Example}
\newtheorem{Definition}{Definition}
\newtheorem{Lemma}[Theorem]{Lemma}
\newtheorem{Proof}{Proof}
\newtheorem{Algorithm}{Algorithm}
\newtheorem{Proposition}{Proposition}
\newtheorem{Corollary}{Corollary}

\begin{abstract} Let   $n, r, e, s$ be are positive integers and the
prime $p,$ the finite local principal ideals ring of parameters
$(p, n, r,e, s)$
$$\texttt{GR}(p^n,r)[x]/(x^e-p\textbf{u}~,~x^s),$$
is defined by an invertible element $\textbf{u}$ of the Galois
ring $\texttt{GR}(p^n,r)$ of characteristic $p^n$ of order
$p^{nr}.$ It is called Galois-Eisenstein ring of parameters $(p,
n, r, e, s).$ A basic problem, which seems to be very difficult is
to determine all non-isomorphism pure Galois-Eisenstein rings of
parameters $(p, n, r, e, s).$ In this paper, this isomorphism
problem for pure Galois-Eisenstein rings of parameters $(p, n, r,
e, s)$ is investigated.
 \end{abstract}

 \vspace{0.15 cm}


\vspace{0.15cm}


\section{Introduction}

Throughout this paper, all rings are finite, associative,
commutative with $1(\neq 0).$ For a ring $R,$ we denote by
$R^\times$ be the set of invertible elements of $R,$ and $J(R)$
the Jacobson radical of $R.$ Let $r, n,$ and $p$ denote positive
integers, $p$ a prime, the residue class ring
$\mathbb{Z}/p^n\mathbb{Z}$ of integers modulo $p^n$ with $p$ prime
and $n > 1$  and let $\texttt{GR}(p^n, r)$ denote the (unique up
to isomorphism) Galois extension of degree $r$ of the ring
$\mathbb{Z}/p^n\mathbb{Z}$ of integers  $\texttt{mod}\,p^n.$  To
begin, let $f\in(\mathbb{Z}/p\mathbb{Z})[x]$ be a primitive
irreducible polynomial of degree $r.$ Then there is a unique monic
polynomial $f_n\in (\mathbb{Z}/p^n\mathbb{Z})[x]$ of degree $r$
such that $f\equiv f_n (\texttt{mod}\,p),$ and $f_n$ divides
$x^{p^r-1} - 1$ in $(\mathbb{Z}/p\mathbb{Z})[x].$ Let $\xi$ be a
root of $f_n,$ so that $\xi^{p^r-1}= 1.$ The Galois ring
$\texttt{GR}(p^n,r)$ is defined to be
$(\mathbb{Z}/p^n\mathbb{Z})[\xi].$ Moreover, $\texttt{GR}(p^n,r)$
is local ring with $J(\texttt{GR}(p^n,r))=(p)$ and the order of
multiplicative subgroup $\Gamma_p(r)^*:=\langle\xi\rangle$ of
\texttt{GR}$(p^n,r)^\times$ is $p^r-1.$ The set
$\Gamma_p(r):=\Gamma_p(r)^*\cup\{0\}$ is called \emph{Teichmüller
set} of \texttt{GR}$(p^n,r).$ The Galois ring $\texttt{GR}(p^n,r)$
depends only on $p,n$ and $r.$

A ring $R$ is a \emph{Galois-Eisenstein ring} of \emph{parameters}
$(p, n, r, e, s)$  if $\texttt{GR}(p^n,r)$ is the largest Galois
ring contained in $R$ and all ideals form the chain \ref{ch1}
\begin{align}\label{ch1}
\{0\}=(\theta^s)\subsetneq(\theta^{s-1})\subsetneq\cdots\subsetneq(\theta)=J(R)\subsetneq
R.
\end{align}
 Since $J(\texttt{GR}(p^n,r))=(p),$ there exists an integer $e$ such that
$(\theta^e)=pR.$ The integer $s$ the \emph{nilpotency index} of
$J(R),$ the Jacobson radical of $R,$ the integer $e$
\emph{ramification index} of $R.$ According to the Theorem 17.5 of
\cite{McDonald}, a GE-ring of parameters $(p, n, r, e, s)$ is
isomorphic to the ring

\begin{align}\label{2}
\texttt{GR}(p^n,r)[x]/(x^e-p\textbf{u}(x)~;~x^s)
\end{align}
 where $\textbf{u}(x):=\textbf{u}_{e-1}x^{e-1}+\cdots+\textbf{u}_1x+\textbf{u}_{0}\in \texttt{GR}(p^n, r)[x],$ with $\textbf{u}_{0}\in \texttt{GR}(p^n, r)^\times.$
 The polynomial $x^e-p\textbf{u}(x)$ is an \emph{associated Eisenstein polynomial}
 to the GE-ring \ref{2}. W. Clark and J. Liang shown in Lemma 2. of \cite{Clark} that when
$n>1,$ an element $\theta$ in $R$ is a root of an Eisenstein
polynomial of degree $e$ over $\texttt{GR}(p^n,r)$ if and only if
$J(R)=(\theta).$ We say that a GE-ring of the form (\ref{2}) is
\emph{pure} if $\textbf{u}(x)=\textbf{u}_0\in
\texttt{GR}(p^n,r)^\times.$

Clark and Liang shown in \cite{Clark} that, if $p\nmid e$ then
pure GE-rings in the form (\ref{2}) are pure and they enumerate
all non-isomorphic pure GE-rings and when $p~|~e,$ Clark and Liang
shown in \cite{Clark} that there are GE-rings which are not pure.
Moreover, Xiang-Dong Hou in \cite{Xiang} gave the number of all
non-isomorphic pure GE-rings, when $n=2$ or $p~|~e,$ $p^2\nmid e$
and $(p-1)\nmid e.$ In this paper, the main goal is the
determination  all non-isomorphic pure GE-ring of parameters
$(p,n,r,e,s),$ when $(p-1)\nmid e.$

 The paper is organized as follows. In Section 2, we
review some basic facts about Galois rings and pure GE-rings of
parameters $(p,n,r,e,s)$ to be used in sequel. In Section 3, we
determine all non-isomorphic pure GE-rings of parameters
$(p,n,r,e,s).$

\section{Preliminaries}

Let $R$ be a pure GE-ring of parameters $(p,n,r,e,s).$ Let
$$\mathcal{L}_e(R):=\left(R^\times\right)^{\cdot
e}\cap\left(\texttt{GR}(p^n,r)\right)^\times$$ be a multiplicative
subgroup of $\texttt{GR}(p^n,r)^\times,$ where
$\left(R^\times\right)^{\cdot e}:=\left\{\textbf{u}^{e}:
\textbf{u}\in R^\times\right\}.$ The aim of this section is the
determination of the integer
$\hyperref[e]{\flat(e)}\in\{1;2;\cdots;n\}$ such that
\begin{align}
\mathcal{L}_e(R)=(\Gamma_p(r)^*)^{\cdot
e}\cap(1+p^{\hyperref[e]{\flat(e)}}\texttt{GR}(p^n,r)).
\end{align}

\subsection{Galois Rings}

The theory of Galois Rings was firstly developed by W. Krull
(1924) and the reader will find in the monograph \cite{Wan}, more
information about on Galois rings quoted. For this subsection, we
gather the results on Galois rings allowing to determine the
integer $\hyperref[e]{\flat(e)}.$ A finite local ring of
characteristic $p^n$ is called Galois ring $\texttt{GR}(p^n,r)$ of
characteristic $p^n$ of rank $r,$ if its Jacobson radical is
generated by $p.$ It is obvious that
$\texttt{GR}(p^n,1)=\mathbb{Z}/p^n\mathbb{Z}$ and
$\texttt{GR}(p,r)=\texttt{GF}(p^r),$ where $\texttt{GF}(p^r)$ is
the Galois field of the size $p^r.$

Let $f_n\in (\mathbb{Z}/p^n\mathbb{Z})[x]$ be a monic polynomial
of degree $r$ such that $$
f_n\;\texttt{mod}\,p\in(\mathbb{Z}/p\mathbb{Z})[x]$$ is
irreducible over $\mathbb{Z}/p\mathbb{Z}$ and $f_n$ divides
$x^{p^r-1} - 1 (\texttt{mod}\, p^n).$ In \cite{Gary}, an algorithm
allows to compute the polynomial $f_n$ is given. Let $\xi$ be a
root of $f_n.$ Then
$\texttt{GR}(p^{n},r)=(\mathbb{Z}/p^n\mathbb{Z})[\xi]$ and
$\Gamma_p(r)^*:=\langle\xi\rangle$ is the unique cyclic subgroup
of order $p^r-1$ of multiplicative group
$\texttt{GR}(p^{n},r)^\times$ isomorphic to cyclic multiplicative
group $\texttt{GF}(p^r)^*.$ The Teichmüller set $\Gamma_p(r)$ of
$\texttt{GR}(p^{n},r)$ forms a complete system of representatives
modulo $p$ in $\texttt{GR}(p^{n},r).$

The following proposition gives the immediate proprieties of the
Teichmüller set of the Galois ring $\texttt{GR}(p^{n},r).$

\begin{Proposition}\label{Lem: Lemma} Let $\xi$ be a generator of $\Gamma_p(r)^*.$ Then

\begin{enumerate}
    \item $\texttt{GR}(p^n,r)\subseteq \texttt{GR}(p^n,r')$ if and
    only if $$\Gamma_p(r)\subseteq\Gamma_p(r')$$ if and
    only if $r$ devises $r';$
    \item $\Gamma_p(r)\cap\Gamma_p(r')=\Gamma_p(gcd(r;r'));$
    \item  $\left(\Gamma_p(r)^*\right)^{\cdot e}=\langle\xi^e\rangle$ and the order of $\langle\xi^e\rangle$ is
$\frac{p^r-1}{gcd(p^r-1;e)}.$
\end{enumerate}
\end{Proposition}

\begin{Example} The monic polynomials $f:=x^2+x+2$ is irreducible over
$\texttt{GF}(3).$ We denote by $\alpha$ the root of $f.$  Then
$\texttt{GF}(3)(\alpha)$ is the Galois field with $9$ elements.
Moreover, the monic polynomial
$f_2:=x^2-5x-1\in(\mathbb{Z}/9\mathbb{Z})[x]$ such
$f_2\;\texttt{mod}=f$ and $f_2$ devises $x^8-1$ in
$\in(\mathbb{Z}/9\mathbb{Z})[x].$ We then construct the Galois
ring $\texttt{GR}(9,2)$ and
$\texttt{GR}(9,2)=(\mathbb{Z}/9\mathbb{Z})[\xi],$ where
$\xi^2=5\xi+1.$ Thus the Teichmüller set of the Galois ring
$\texttt{GR}(9,2)$ is
$$\Gamma_3(2)=\{0,1,\xi,\xi^2(=5\xi+1),\xi^3(=8\xi+5),\xi^4(=8),\xi^5(=8\xi),\xi^6(=4\xi+8),\xi^7(=\xi+4)\},$$

and $\left(\Gamma_3(2)^*\right)^{\cdot
18}=\left(\Gamma_3(2)^*\right)^{\cdot
2}=\langle\xi^2\rangle=\{1,5\xi+1,8\}.$
\end{Example}

 The following lemma gives the immediate proprieties of a complete residue system modulo
 $p^{\ell},$ where $\ell\in\{0,1,\cdots,n\}.$

\begin{Lemma}\label{Lem: Lemma} Let $\xi$ be a generator of $\Gamma_p(r)^*$ and $\ell\in\{0,1,\cdots,n\},$  we
consider the set
\begin{align}\mathcal{R}_r(\ell):=\left\{\sum_{i=0}^{\ell-1}\xi_ip^i\in
\texttt{GR}(p^{n},r) :\xi_i\in\Gamma_p(r)\right\}\end{align} and
by convention, we adopt $\mathcal{R}_r(0):=\{0\}.$ Then

\begin{enumerate}
    \item $\mathcal{R}_r(1)=\Gamma_p(r)$ and $\mathcal{R}_r(n)=\texttt{GR}(p^{n},r);$
    \item $\mathcal{R}_r(\ell)$ forms a complete residue system modulo $p^{\ell}$ in $\texttt{GR}(p^{n},r);$
    \item \label{item4} for each  $\alpha\in\texttt{GR}(p^{n},r),$ for each  $\ell\in\{0;1;2;\cdots;n\},$ there
    exists a unique $(\gamma;\beta)$ in $\mathcal{R}_r(\ell)\times
\mathcal{R}_r(n-\ell),$ such that $\alpha=\gamma+p^\ell\beta.$

    \end{enumerate}
\end{Lemma}

We remark that
\begin{align}
\mathcal{R}_r(0)\subset\mathcal{R}_r(1)\subset\cdots\subset\mathcal{R}_r(n).
\end{align}

\begin{Proposition} The automorphism group of the ring
$\texttt{GR}(p^{n},r)$ is
$$\texttt{Aut}(\texttt{GR}(p^{n},r)):=\langle\{\sigma_p :\xi\mapsto
\xi^p\}\rangle,$$ and for all $\xi_i\in\Gamma_p(r),$

$$\sigma\left(\sum_{i=0}^{n-1}\xi_ip^i\right)=\sum_{i=0}^{n-1}\sigma(\xi_i)p^i.$$

 Moreover, the groups
$\texttt{Aut}(\texttt{GR}(p^{n},r))$ and $\{0; 1;\cdots;r-1\},$
are isomorphic and for all $j\in\{0; 1; \cdots; r-1\},$
$$\left\{x\in\Gamma_p(r):
\sigma_p^j(x)=x\right\}=\Gamma_p\left(gcd(r,j)\right).$$
 \end{Proposition}
We remark that if $gcd(r,j)=1$ then $\left\{x\in\Gamma_p(r):
\sigma_p^j(x)=x\right\}=\mathbb{Z}/p^n\mathbb{Z}.$

\begin{Example} Consider the Galois
ring $$\texttt{GR}(9,2)=(\mathbb{Z}/9\mathbb{Z})[\xi],$$ where
$\xi^2=5\xi+1.$ Thus the Teichmüller set of the Galois ring
$\texttt{GR}(9,2)$ is
$$\Gamma_3(2)=\{0,1,\xi,5\xi+1,8\xi+5,8,8\xi,4\xi+8,\xi+4\}.$$

The Galois group of of the ring $\texttt{GR}(p^{n},r)$ is
$$\texttt{Aut}(\texttt{GR}(9,2))=\{Id, \sigma_3:\xi\mapsto
\xi^3\}$$

and

$$\left\{x\in\Gamma_3(2):
\sigma_3(x)=x\right\}=\Gamma_3\left(gcd(2,1)\right)=\mathbb{Z}/9\mathbb{Z}.$$
\end{Example}

The structure of group of invertible elements of
$\texttt{GR}(p^{n},r)$ is given by the following theorem.

\begin{Theorem}([\cite{Wan} Theorem 14.11])\label{thm2.2} Let $\texttt{GR}(p^{n},r)^\times$ be the
group of invertible elements of $\texttt{GR}(p^{n},r).$ Then
$\texttt{GR}(p^{n},r)^\times$ is the internal direct product of
subgroups $\Gamma_p(r)^*$ and $1+p\texttt{GR}(p^{n},r).$ Moreover,
\begin{enumerate}
    \item If $p$ is odd or if $p=2$ and $n\leq 2,$ then
    \begin{align}1+p\texttt{GR}(p^{n},r)\cong
    \left(\mathbb{Z}/(p^{n})\right)^r,\end{align}
    \item if $p=2$ and $n\geq 3,$ then
    \begin{align}1+p\texttt{GR}(p^{n},r)\cong
    \texttt{GF}(2)\times\left(\mathbb{Z}/(2^{n-2})\right)\times\left(\mathbb{Z}/(2^{n-1})\right)^{r-1}.\end{align}
\end{enumerate}
\end{Theorem}

The Theorem \ref{thm2.2} has as a consequence the following
corollary and this corollary gives a simple expression of the
subgroup $(1+p\texttt{GR}(p^{n},r))^{\cdot p^i}$ of
$1+p\texttt{GR}(p^{n},r).$

\begin{Corollary}
Let $i$ be an integer and  $p$ is a prime. If $p$ odd or if $p=2$
and $n\leq 2,$ then
\begin{align}\label{cor}(1+p\texttt{GR}(p^{n},r))^{\cdot
p^i}=1+p^{i+1}\texttt{GR}(p^{n},r).
\end{align}
\end{Corollary}

\subsection{Pure Galois-Eisenstein Rings}

Let $R$ be a pure Galois-Eisenstein ring (short: pure GE-ring) of
parameters $(p, n, r, e, s).$ Then there exists $\textbf{u}\in
\texttt{GR}(p^n,r)^\times$ such that
$R=\texttt{GR}(p^n,r)[x]/(x^e-p\textbf{u}~;~x^s).$ The writing
\begin{align}
\texttt{GE}(\textbf{u}):=\texttt{GR}(p^n,r)[x]/(x^e-p\textbf{u}~;~x^s)
\end{align}
means that $\texttt{GE}(\textbf{u})$ is a GE-ring of parameters
$(p,n,r,e,s)$ defines by the invertible element
$\textbf{u}\in\texttt{GR}(p^n,r)^\times.$ In the sequel, we write
$\texttt{GE}(\textbf{u})=\texttt{GR}(p^n,r)[\theta],$ such that
$\theta^e=p\textbf{u}$ and $\theta^s=0,$ but $\theta^{s-1}\neq 0,$
for denote the pure Galois-Eisenstein ring of parameters $(p, n,
r, e, s).$ We obvious that the Jacobson radical of
$\texttt{GE}(\textbf{u})$ is $(\theta).$

\begin{Lemma}\label{lem2.3} Let $R$ be a pure Galois-Eisenstein ring of parameters $(p, n, r,
e, s).$ Then the group of invertible elements $R^\times$ of $R$ is
the internal direct product of subgroups $\Gamma_p(r)^*$ and
$1+J(R).$
\end{Lemma}

The following theorem gives the structure of subset
\begin{align}
\mathcal{L}_e(R):=\left(R^\times\right)^{\cdot
e}\cap\left(\texttt{GR}(p^{n},r)\right)^\times
\end{align}
 of group of invertible elements
$\texttt{GR}(p^{n},r)^\times.$  This is the main result in this
section.

\begin{Theorem}\label{Thm2} Let $R$ be a pure Galois-Eisenstein ring of parameters $(p, n, r,
e, s)$ such that $s=(n-1)e+t,$ $0\leq t \leq e$ and $(p-1)\nmid
e.$ Then there exists an integer
$\label{e}\hyperref[e]{\flat(e)}\in\{1;2;\cdots;n\}$ such that
$$\mathcal{L}_e(R)=\left\langle\eta\right\rangle\cdot(1+p^{\hyperref[e]{\flat(e)}}\texttt{GR}(p^{n},r)),$$
where
$$\hyperref[e]{\hyperref[e]{\flat(e)}}= \left\{%
\begin{array}{ll}
    1, & \hbox{if $t\leq \frac{e}{p}$ and $n=2$;} \\
    min\{\omega+1;n\}, & \hbox{if $t>\frac{e}{p}$ and $n= 2$ or $n\neq 2,$}, \\
\end{array}
\right.$$ where $\omega:=max\{i\in\mathbb{N}^*:p^i~|~e\}.$
\end{Theorem}

\begin{Proof}
The size of the multiplicative group $1+J(R)$ is a power of $p,$
therefore

$$(1+J(R))^{\cdot
p^{e}}=1+J(R))^{\cdot p^{\omega}}$$ and
$\hyperref[e]{\flat(e)}=\hyperref[e]{\flat(\omega)}.$

By \hyperref[thm2.2;em2.3]{Theorem 2.2 and Lemma 2.3.},
$$\mathcal{L}_e(R)=\left((\Gamma_p(r)\cdot(1+J(R))^{\cdot
p^{\omega}}\right)\cap\left(\Gamma_p(r)\cdot(1+p\texttt{GR}(p^{n},r)\right),$$
It follows that
$$\mathcal{L}_e(R)=\left\langle\eta\right\rangle\cdot\left((1+J(R))^{\cdot
p^{\omega}}\cap(1+p\texttt{GR}(p^{n},r))\right),$$ in according by
\hyperref[Lem: Lemma]{Lemma 2.1.}, we have
$\left(\Gamma_p(r)^*\right)^{\cdot
e}=\left\langle\eta\right\rangle.$  It suffices to determine the
integer $\hyperref[e]{\flat(e)}$ such that

$$(1+J(R))^{\cdot
p^{\omega}}\cap(1+p\texttt{GR}(p^{n},r))=1+p^{\hyperref[e]{\flat(e)}}\texttt{GR}(p^{n},r)).$$

We write

$$(1+J(R))^{\cdot
p^\omega}\cap(1+p\texttt{GR}(p^{n},r))=\mathcal{L}_1\cup
\mathcal{L}_2,$$

 where

$$ \mathcal{L}_1=\left\{(1+\theta^b\varepsilon)^{p^\omega}\in
1+p\texttt{GR}(p^{n},r)~:~\geq e~~;~~ \varepsilon\in
R^\times\right\},$$
$$\mathcal{L}_2=\left\{(1+\theta^b\varepsilon)^{p^\omega}\in 1+p\texttt{GR}(p^{n},r)~:~ 1\leq b<
e~~;~~  \varepsilon\in R^\times\right\}.$$

Let $b$ be an integer and $\varepsilon\in R^\times.$ On the one
hand, suppose that $b\geq e.$ We have

 $$1+\theta^b\varepsilon\in 1+pR$$

and

$$(1+\theta^b\varepsilon)^{p^\omega}\in (1+pR)^{\cdot p^\omega}\subseteq 1+p^{\omega+1}R.$$

 Thus

\begin{eqnarray*}
 (1+\theta^b\varepsilon)^{p^\omega}\in(1+p^{\omega+1}R)\cap
(1+p\texttt{GR}(p^{n},r))&=&1+p^{\omega+1}\texttt{GR}(p^{n},r) \\
   &=& (1+p\texttt{GR}(p^{n},r))^{\cdot
p^\omega}.
\end{eqnarray*}

Since

$$(1+p\texttt{GR}(p^{n},r))^{\cdot p^\omega}\subseteq
\mathcal{L}_1,$$

 we have

$$\mathcal{L}_1=(1+p\texttt{GR}(p^{n},r))^{\cdot p^\omega},$$

 on the
other hand, suppose that $b<e,$ we develop the following
expression $(1+\theta^b\varepsilon)^{p^\omega}$ and we obtain:
\begin{eqnarray*}
  (1+\theta^b\varepsilon)^{p^\omega} &=& 1+\sum_{l=1}^{p^\omega-1}\left(\begin{array}{c} l \\  p^\omega \\\end{array}\right)(\theta^b\varepsilon)^l+ (\theta^b\varepsilon)^{p^\omega},\text{since $l=jp^i$ and $p~\nmid~j$}; \\
    &=&  1+ \sum_{i=1}^{\omega-1}\sum\limits_{\substack{  jp^i\leq p^\omega-1 \\ p~\nmid~j}}\left(\begin{array}{c} jp^i\\p^\omega\\\end{array}\right)(\theta^b\varepsilon)^{jp^i}+(\theta^b\varepsilon)^{p^\omega},\\
\end{eqnarray*}

      since

     $$p^{\omega-i}|\left(\begin{array}{c} jp^i\\ p^\omega\\\end{array}\right)$$

     and

     $$p^{\omega-i+1}\nmid\left(\begin{array}{c} jp^i\\
      p^\omega\\\end{array}\right);$$

 thus

$$(1+\theta^b\varepsilon)^{p^\omega} = 1+
\sum_{i=1}^{\omega-1}\theta^{e(\omega-i)+bp^i}\varepsilon_i +
\theta^{bp^\omega}\varepsilon^{p^\omega},$$

 where $\varepsilon_{i}\in R^\times$. We write $h(b,i):=e(\omega-i)+bp^i$ and
$$h(b):=min\{e(\omega-i)+bp^i:0\leq i< \omega\}.$$ If $(p-1)\nmid e,$ then $h(b,i)\neq h(b,i+1).$ Thus there exists a unique
$\imath\in\{0;\cdots;\omega-1 \}$ such that $h(b)=h(b,\imath).$
Therefore, there exists $\varepsilon_b\in R^\times$ such that
$$(1+\theta^b\varepsilon)^{p^\omega}=1+ \theta^{h(b)}\varepsilon_b
+ \theta^{bp^\omega}\varepsilon^{p^\omega}.$$

This forces that $h(b)\geq s,$ and $bp^\omega\equiv 0
~\texttt{mod}~ e.$  We have $a\in \mathbb{N}^*$ such that $
(p\textbf{u})^a\varepsilon^{p^\omega}\in
 p\texttt{GR}(p^{n},r).$ We obtain,
$$\mathcal{L}_2=
\left\{%
\begin{array}{ll}
    1+p\texttt{GR}(p^{n},r), & \hbox{if $n=2$  and $t\leq\frac{e}{p}$;} \\
    \{1\}, & \hbox{if $n\neq 2$ or $n=2$ and $t> \frac{e}{p}.$}~~~~\blacksquare

\end{array}%
\right. $$

\end{Proof}

\begin{Example} Let $R$ be a GE-ring of parameters $(3,2,2,18,25).$

Then by \hyperref[Thm2]{Theorem 2.2}, $\flat(2)=2.$ Thus
\begin{eqnarray*}
  \mathcal{L}_{9}(R) &=& \langle\xi^2\rangle\cdot(1+3^8\texttt{GR}(9,2)) \\
      &=& \{1,\xi^2,\xi^4\},
\end{eqnarray*}
where $\Gamma_3(2)^*=\langle\xi\rangle$ and $\xi^2 = 5\xi + 1.$

\end{Example}

\section{Main result}

In this section, we determine the pure GE-rings of parameters
$(p,n,r,e,s).$ Note when $n=1,$ $R=GF(p^r)/(x^e).$ Such rings need
no classification, so we will suppose that $n\geq 2.$ The
isomorphism problem for pure GE-rings with given parameters
$(p,n,r,e,s)$ mentions by A. A. Nechaev in \cite{Nechaev} and T.
G. Gazaryan in \cite{T. G. Gazaryan}, shown that pure GE-rings
\texttt{GE}$(1)$ and \texttt{GE}$(2)$ of parameters $(3,2,1,2,3)$
are non-isomorphic rings with isomorphic additive and
multiplicative groups.

 Let $\sim$ be the equivalence relation defined on $\texttt{GR}(p^n,r)^\times$ by

\begin{align}
 \textbf{u} \sim \textbf{v}  \Leftrightarrow \texttt{GE}(\textbf{u})\cong \texttt{GE}(\textbf{v}),
\end{align}
where $\textbf{u},\textbf{v}\in \texttt{GR}(p^n,r)^\times.$

Let $\sigma$ be a generator of
$\texttt{Aut}(\texttt{GR}(p^{n},r)),$ and the multiplicative
subgroup
$$\mathcal{L}_e(R)=\left\langle\eta\right\rangle\cdot(1+p^{\hyperref[e]{\flat(e)}}\texttt{GR}(p^{n},r))$$
of $\texttt{GR}(p^n,r)^\times$ and $\eta$ a generator of
$\left(\Gamma_p(r)^*\right)^{\cdot e}.$ Consider
\begin{align}
\texttt{U}_r(\hyperref[e]{\flat(e)}):=\langle{\chi}\rangle+p\mathcal{R}_r(\hyperref[d]{\partial(e)}-1),
\end{align}
where $\partial(e)=min\{\partial(e);n-1\}.$

We write $\chi:=\xi^{\frac{p^r-1}{d}}$ and $d:=gcd(p^r-1;e).$
Since $gcd\left(d;(\frac{p^{r}-1}{d})\right)= 1,$ then
$\langle\chi\rangle\cdot\langle\eta\rangle=\Gamma_p(r)^*.$

The isomorphic GE-rings have the same parameters, but there are
non-isomorphic GE-ring with the same parameters. The following
example illustrates the isomorphism problem for the GE-rings.

\begin{Example} T. G. Gazaryan, given in \cite{T. G. Gazaryan}, two GE-rings of
parameters $(3,2,1,2,3)$
$\texttt{GE}(1)=(\mathbb{Z}/9\mathbb{Z})[\theta]$ and
$\texttt{GE}(2)=(\mathbb{Z}/9\mathbb{Z})[\delta]$  as an example
of non-isomorphic commutative chain rings. We have $\theta^2=3,
\delta^2=6,$ and $\theta^3=\delta^3=0.$

  Then the radical Jacobson are:
$$J(\texttt{GE}(1))=(\theta)=\{0;3;6;\theta;3\theta\}$$ and
$$J(\texttt{GE}(1))=(\delta)=\{0;3;6;\delta;3\delta\}.$$
The Teichmüller set $\Gamma_3(1)=\{0;1;8\},$ allows to write
$$\texttt{GE}(1)=\{a+b\theta: (a;b)\in\Gamma_3(1)\}$$ and
$$\texttt{GE}(2)=\{a+b\delta: (a;b)\in\Gamma_3(1)\}.$$ But
$J(\texttt{GE}(1))\neq J(\texttt{GE}(2)).$ Indeed, if
$J(\texttt{GE}(1))= J(\texttt{GE}(2)),$ then there exists $z\in
\texttt{GE}(1)^\times$ such that $\theta=z\delta.$ As
$$3=\theta^2=z^2\delta^2=6z^2.$$ It follows, $z^2=(a+\theta
b)^2=a^2+2\theta ab.$ Therefore, the equation $3=6z^2$ is
equivalent to $3=6a^2.$ Now, for all $a\in
(\mathbb{Z}/9\mathbb{Z})^\times=\{1;2;4;5;7;8\},$ we have
$6a^2=6.$ It is absurd.
\end{Example}
\begin{Lemma}\label{Lem'} Let $R$ be the GE-ring of
parameters $(p,n,r,e,s).$ Then there exists $\textbf{v}$ in
$\texttt{U}_r(\hyperref[e]{\flat(e)})$ such that
$R=\texttt{GE}(\textbf{v}).$
\end{Lemma}

\begin{Proof} Let $R$ be the GE-ring of
parameters $(p,n,r,e,s).$ Then there exists $\textbf{u}\in
\texttt{GR}(p^n,r)^\times$ such that $R=\texttt{GE}(\textbf{u}).$
By the item\ref{item4} of \hyperref[Lem: Lemma]{Lemma 2.1.}, there
exists
$$(\gamma;\beta)\in\mathcal{R}_r(\hyperref[e]{\flat(e)})\times
\mathcal{R}_r(n-\hyperref[e]{\flat(e)}),$$ such that
$\textbf{u}=\gamma+p^{\hyperref[e]{\flat(e)}}\beta$ and $\gamma\in
\texttt{GR}(p^n,r)^\times.$ Since
$\langle\delta\rangle\cdot\langle\eta\rangle=\Gamma_p(r)^*,$ there
exists $(\alpha,\textbf{v})\in\langle\eta\rangle\times
\texttt{U}_r(\hyperref[e]{\flat(e)})$ such that $\gamma=\alpha
\textbf{v}.$

Now, $\textbf{u}=
\textbf{v}(\alpha+p^{\hyperref[e]{\flat(e)}}\beta),$ where
$\beta=\beta\textbf{v}^{-1}.$ We have
$$g:=\alpha+p^{\hyperref[e]{\flat(e)}}\beta\in \mathcal{L}_e(R).$$
By \hyperref[Thm2]{Theorem 2.4,} there exists $z\in R^\times$ such
that $g =z^e.$

Since

$$\texttt{GE}(\textbf{u})=\texttt{GR}(p^n,r)[\theta]$$ with
$\theta^e= p\textbf{v}$  and $s=min\{i\in\mathbb{N}:\theta^i=0\}.$
Thus, $J(\texttt{GE}(\textbf{u}))=(\delta)$ where
$\delta:=z\theta.$

We have $\delta^e=p\textbf{v}$ and
$s=min\{i\in\mathbb{N}:\delta^i=0\}.$ So
$\texttt{GE}(\textbf{u})=\texttt{GE}(\textbf{v}).~\blacksquare$
\end{Proof}

 The following
lemma gives a fundamental condition for non-isomorphic pure
GE-rings.

\begin{Lemma}\label{Lem} For each $\textbf{u}$ and $\textbf{v}$ in
$\texttt{U}_r(\hyperref[e]{\flat(e)}).$ Then
$$ \textbf{u} \sim \textbf{v} \Leftrightarrow
\textbf{u}=\sigma^i(\textbf{v}),$$ for some $i\in \{0;
1;\cdots;r-1\}.$
\end{Lemma}

\begin{Proof} Let $\theta:=x+(x^e-p\textbf{u}~;~x^s)$ be a generator of
$J(\texttt{GE}(\textbf{u})).$  Then pure GE-rings
$\texttt{GE}(\textbf{u})$ and $\texttt{GE}(\textbf{v})$ are
isomorphic means by Lemma. 4 of \cite{Clark}, the existence of $
(i; z)\in \{0; 1;\cdots;r-1\}\times \texttt{GE}(\textbf{u})^\times
$ such that $(\theta z)^e=p\sigma^i(\textbf{v})$ and
$\theta^e=p\textbf{u}.$ Since $\textbf{u}z^e=\sigma^i(\textbf{v})$
and $\textbf{u},
\sigma^i(\textbf{v})\in\texttt{U}_r(\hyperref[e]{\flat(e)}),$ we
have $z^e\in1+p^{\hyperref[e]{\flat(e)}}\texttt{GR}(p^{n},r). $
Therefore,
\begin{eqnarray*}
 \textbf{u} \sim \textbf{v}  &\Leftrightarrow & \exists (i;\gamma)\in \{0; 1;\cdots;r-1\}\times(1+p^{\hyperref[e]{\flat(e)}}\texttt{GR}(p^{n},r)) : \textbf{u}\gamma=\sigma^i(\textbf{v}),~\text{where}~\gamma:=z^e; \\
    &\Leftrightarrow&  \exists (i;\gamma)\in \{0; 1;\cdots;r-1\}\times (1+p^{\hyperref[e]{\flat(e)}}\texttt{GR}(p^{n},r)):
    \textbf{u}\gamma=\sigma^i(\textbf{v}).
\end{eqnarray*}
Now, $\gamma\in
(1+p^{\hyperref[e]{\flat(e)}}\texttt{GR}(p^{n},r))\cap\texttt{U}_r(\hyperref[e]{\flat(e)})=\{1\}.$
$\blacksquare$
\end{Proof}

Let $\textbf{u}\in\texttt{U}_r(\hyperref[e]{\flat(e)}),$ we write
$\texttt{C}(\textbf{u})=\{\sigma^i(\textbf{u}): i\in \{0;
1;\cdots;r-1\}\},$ the Frobenius class of $\textbf{u}$ and
$\mathcal{C}_r(\hyperref[e]{\flat(e)})$ a complete set of
Frobenius representative of
$\texttt{U}_r(\hyperref[e]{\flat(e)}).$

\begin{Theorem}\label{thm}  Let $\Xi^*(p, n, r, e, s)$ be the set of pure GE-rings of parameters
$(p, n, r, e, s)$ up to isomorphism. Suppose that $(p-1)\nmid e.$
Then the mapping
\begin{align}\begin{array}{ccccc} \texttt{GE} : &\mathcal{C}_r(\hyperref[e]{\flat(e)}) &\to &  \Xi^*(p, n, r, e, s)\\
& \textbf{u} & \mapsto & \texttt{GE}(\textbf{u})
\end{array}\end{align} is bijective.
 Moreover,
\begin{align}
|\Xi^*(p, n, r, e,
s)|=\frac{1}{r}\sum_{i=1}^{r}gcd(p^{gcd(r,i)}-1,e)p^{({\hyperref[d]{\partial(e)}-1})gcd(i,r)}.
\end{align}
\end{Theorem}

\begin{Proof} By \hyperref[Lem']{Lemma 3.1}, the mapping $\texttt{GE}$ is bijection.
Indeed, For each
$\textbf{u}\in\texttt{U}_r(\hyperref[e]{\flat(e)}),$ there exists
a unique $\textbf{v}$ in $\mathcal{C}_r(\hyperref[e]{\flat(e)})$
such that $\texttt{GE}(\textbf{u})=\texttt{GE}(\textbf{v}).$

 Now, we determine the size of\,  $\Xi^*(p, n, r, e, s).$
Then we use the action of $\sigma$ on
$\texttt{U}_r(\hyperref[e]{\flat(e)})$ is given by:
\begin{align}
\begin{array}{ccccc} \sigma: &\texttt{U}_r(\hyperref[e]{\flat(e)}) & \to &\texttt{U}_r(\hyperref[e]{\flat(e)}) \\
& \varepsilon & \mapsto &\sigma(\varepsilon).
\end{array}
\end{align}
 For each $1\leq i \leq r,$ we have:
 \begin{eqnarray*}
   fix(\sigma^i)&=&\left\{\varepsilon\in \texttt{U}_r(\hyperref[e]{\flat(e)}): \sigma^i(\varepsilon)=\varepsilon \right\}; \\
      &=&
      \texttt{U}_{r_i}(\hyperref[e]{\flat(e)}), \; \text{where}\; r_i:=gcd(r,i).
 \end{eqnarray*}
Thus by the Burnside's Lemma, the number of
$\langle\sigma\rangle-$orbits in
$\texttt{U}_r(\hyperref[e]{\flat(e)})$ is
\begin{align*}
  \frac{1}{|\langle\sigma\rangle|}\sum_{i=1}^{r}|fix(\sigma^i)|=\frac{1}{r}\sum_{i=1}^{r}gcd(p^{r_i}-1,e)p^{({\hyperref[d]{\partial(e)}-1})r_i}. ~~~~~~~~~~\blacksquare
\end{align*}
\end{Proof}

Thus, the mapping $\texttt{GE}$ allows to determine up to
isomorphism  all pure GE-rings.

\begin{Example} Consider the pure GE-rings
of parameters $(3,2,1,2,3).$ In according by Theorem
\hyperref[Thm2]{Theorem 2.3,} $\flat(2)= 1.$  Then
$$\mathcal{C}_1(1)=\texttt{U}_1(1)=\{1; 8\}.$$
By formula of \hyperref[Thm]{Theorem 3.3,} there are two
non-isomorphism pure GE-rings of parameters $(3,2,1,2,3).$
Therefore, in the article \cite{T. G. Gazaryan}, the GE-rings
$r(1)$ et $r(8)$ are only non-isomorphism pure GE-rings of
parameters $(3,2,1,2,3).$
\end{Example}

\begin{Example} Consider the GE-rings $\texttt{GR}(9,2)[\theta]$ of parameters
$(3,2,2,18,25)$ and
$\texttt{GR}(9,2)=(\mathbb{Z}/9\mathbb{Z})[\xi],$ with
$\xi^2=5\xi+1.$

 In according by Theorem \hyperref[Thm2]{Theorem
2.3,} $\flat(18)= 2.$ Then
$$\mathcal{C}_2(18)=\texttt{U}_2(2)=\langle\xi^4\rangle=\{1;8\}.$$
 By
formula of \hyperref[Thm]{Theorem 3.3,} there are 2
non-isomorphism pure GE-rings of parameters $(3,2,2,18,25),$ and
the GE-rings of parameters $(3,2,2,18,25)$ up to isomorphism are
$r(1),$ and $r(8)$
\end{Example}

\section*{Conclusion}

A pure Galois-Eisenstein ring of $(p,n,r,e,s)$ is constructed from
an element of the complete set of Frobenius representative of
$\texttt{U}_r(\hyperref[e]{\flat(e)})$ and this construction are
unique up to isomorphism. In general, the isomorphism problem for
Galois-Eisenstein Rings stays open.

\end{document}